\documentclass[11pt]{amsart}
\usepackage{graphics}
\usepackage{amsmath,amsthm,amscd}
\usepackage{amssymb,amsmath}
\usepackage{amscd}
\usepackage{epsfig}
\usepackage[all,cmtip]{xy}
\usepackage{xypic}
\newcommand{\Ext}{\mathrm{Ext}}
\newcommand{\Hom}{\mathrm{Hom}}

\newcommand{\rad}{\mathrm{rad}}

\newcommand{\e}{\varepsilon}
\newcommand{\cal}{\mathcal}

\newcommand{\Ker}{\mathrm{Ker}}
\newcommand{\Ima}{\mathrm{Im}}
\newcommand{\HH}{\mathrm{HH}}

\def\aa{{\mathfrak a}}
\def\bb{{\mathfrak b}}
\def\cc{{\mathfrak c}}
\def\dd{{\mathfrak d}}

\newtheorem{Theorem}{Theorem}[section]
\newtheorem{Lemma}[Theorem]{Lemma}%[section]
\newtheorem{Proposition}[Theorem]{Proposition}%[section]
\newtheorem{Definition}[Theorem]{Definition}%[section]

\begin{document}

\title{On the Hochschild cohomology of tame Hecke algebras }

\thanks{The second author acknowledges support through an EPSRC Postdoctoral Fellowship EP/D077656/1 as well as through a Leverhulme Early Career Fellowship}
\author{Karin Erdmann and Sibylle Schroll}
\address{Mathematical Institute, University of Oxford, 24-29 St. Giles',
Oxford OX1 3LB, UK\newline Department of Mathematics, University of
Leicester, University Road, Leicester LE1 7RH, UK }
\email{erdmann@maths.ox.ac.uk, schroll@mcs.le.ac.uk}
\begin{abstract}
We explicitly calculate a projective bimodule resolution for a
special biserial algebra giving rise to the Hecke algebra ${\mathcal
H}_q(S_4)$ when $q=-1$. We then determine the dimensions of the
Hochschild cohomology groups.
\end{abstract}
\maketitle

\hsize440pt\hoffset-30pt
 \baselineskip=16pt

 \section{Introduction}

\parindent0pt

In this paper we are interested in Hochschild cohomology of
finite-dimensional algebras; the main motivation is to generalize
group cohomology to larger classes of algebras. If suitable finite
generation holds, one can define support varieties of modules as
introduced by \cite{SS}. Furthermore, when the algebra is
self-injective, many of the properties of group representations
generalize to this setting as was shown in \cite{EHSST}. Although
Hecke algebras do not have a Hopf algebra structure, one may expect
that being deformations of group algebras,  they should have good
homological properties. Furthermore, the study of Hochschild
cohomology for blocks of group algebras with cyclic or dihedral
defect groups has made use of the fact that the basic algebras in
this case are special biserial \cite{Ho1,Ho2}. These results suggest
that Hochschild cohomology of special biserial algebras might be
accessible more generally.
This paper should be seen as a step towards understanding Hochschild
cohomology of special biserial algebras.

\ \
Here we deal with the self-injective special biserial
algebra $A$ which occurs as the basic algebra of the Hecke algebra
${\cal H}_q(S_4)$ when $q=-1$ (and the characteristic of the field is $\neq 2$). For this algebra we will explicitly
give a minimal projective bimodule resolution and use it to
calculate the dimensions of the  Hochschild cohomology. By
\cite{Erdmann-Nakano} there are two Morita types of tame blocks for  Hecke algebras ${\cal
H}_q(S_n)$, and the other one   is derived equivalent to $A$;
this can be seen directly, and it also follows from \cite{CR}. Since Hochschild cohomology
is invariant under derived equivalence, our result gives information
for arbitrary tame blocks of Hecke algebras of type $A$.

\ \
In order for the Hochschild cohomology ring to be finitely
generated, the dimensions of $HH^n(A)$ must have at most polynomial
growth, and this follows from our result.
In \cite{Erdmann-Solberg}, it is shown by a different
approach, via derived equivalence, that the Hochschild cohomology ring
in this case must be finitely generated and that
the finite generation hypothesis in \cite{EHSST} holds, but
the approach does not give insight in the structure of $HH^*(A)$.

\ \ The algebra $A$ can be described by a quiver and relations (cf.
\cite{Erdmann-Nakano}). Let $K$ be any  field and let ${\cal Q}$ be
the quiver
$$\xymatrix{
1\ar@(ul,dl)[]_{\varepsilon}\ar@<1ex>[r]^{\alpha} & 2
\ar@<1ex>[l]^{\beta} \\ }$$
%\begin{figure}[h]\label{Quiver}
%\begin{center}
% \includegraphics[ width=4.5cm]{Quiver.eps}
%\end{center}
%\caption{\small The quiver of $A$}
%\end{figure}
with relations
\begin{equation} \e \alpha = \beta \e = 0, (\alpha \beta)^2 = \e^2.  \hfill \end{equation}
Then $A=K{\cal Q}/I$ where $I$ is generated by the relations in (1).
We denote the trivial path at the vertex $i$, $i=1,2$, by $e_i$ and
the corresponding simple $A$-module by $S_i$. We read paths from
left to right and, unless stated differently, all the modules
considered will
 be (finitely generated) right $A$-modules. Let
$\rad A$ be the Jacobson radical of $A$ and denote the enveloping
algebra of $A$ by $A^e = A^{op} \otimes_k A$. For a path $a$ in $A$
we denote by ${\bf o}(a)$ the idempotent at the vertex from which
the path leaves and by ${\bf t}(a)$ the idempotent at the vertex at
which the path arrives. We say an element $p$ in $K{\mathcal Q}$ is
right (resp. left) uniform if it is a linear combination of paths
ending in (resp. beginning at) a  single vertex. An element $p$ in
$K{\mathcal Q}$ is uniform if it is left and right uniform.

\bigskip
\section{Construction of a minimal projective bimodule resolution}
 In this section we construct a minimal projective $A$-$A$-bimodule
resolution of $A$, that is, we construct an exact complex
$$R_\bullet: \;\;\;\;\;\;\;  \ldots \longrightarrow R_n \stackrel{\delta_n}{\longrightarrow} R_{n-1} \stackrel{\delta_{n-1}}{\longrightarrow} \ldots
\longrightarrow R_1 \stackrel{\delta_1}{\longrightarrow} R_0
\stackrel{\delta_0}{\longrightarrow} A \longrightarrow 0$$ with $R_n
= \bigoplus P_{ij}$ and where $P_{ij}$ is the projective
$A$-$A$-bimodule $Ae_i \otimes e_j A$.  By \cite{Happel} the projective
 $P_{ij}$  occurs in $R_n$ exactly as many times as ${\rm dim}
\Ext^n_A(S_i, S_j)$. Recall that $\Ext^n_A(S_i, S_j)\cong
\Hom_A(\Omega^n(S_i), S_j)$ where $\Omega(X)$ denotes the kernel of
a projective cover of $X$. These dimensions can be found in at least
two ways. First, one can use the method as outlined in the appendix
of \cite{R} to calculate explicitly the $\Omega$-translates of the
simple modules, and read off the dimensions of ${\rm Ext}^n_A(S_i,
S_j)$. Second, they also follow from the construction of a mininimal
projective resolution, following the method of
\cite{Green-Solberg-Zacharia}, which we give in 2.2. Namely, the
dimension of $\Ext^n_A(S_i, S_j)$ is the number of elements $x$ in
$g_n$ satisfying $x=e_ixe_j$.

The answer is  as follows. Write $n= 4k+l$ where  $k \geq 0$ and $0 \leq l
\leq 3$, then

$$\dim \Ext^{4k+l}_A(S_1, S_1) =
\left\{ \begin{array}{ll}
 2k+1 & \mbox{if $ l = 0,1,2$} \\
 2k+2 & \mbox{if $ l = 3$} \\
 \end{array} \right.        $$
and
$$\dim \Ext^{4k+l}_A(S_i, S_j) =
\left\{ \begin{array}{ll}
 0 & \;\;\;\; \mbox{if $ l = 0,3$ and  $i \neq j$ \;\; or \;\; if $l = 1,2$ and $i = j =2$} \\
 1 & \;\;\;\; \mbox{if $ l = 1,2$ and $i \neq j$ \;\; or \;\; if $ l = 0,3$ and $i = j =2$. } \\
 \end{array} \right.$$

%Minimal projective resolutions of simple modules
%in this case can be calculated, for example similar to the Appendix in
%\cite{R}.
%For the algebra $A$, the dimensions of the projective modules in a
%minimal projective resolution of a simple module are of linear
%growth. Combined with the fact that the algebra only has two
%indecomposable projectives this easily gives formulas of growth for
%the dimensions of the $\Ext^n(S_i, S_j)$.
Then by Happel's result the
terms of the projective bimodule resolution $(R_\bullet, \delta)$
are  for $k\geq 0$,
\begin{equation}\label{terms}\begin{array}{rl} R_{4k} \cong & P_{11}^{2k+1} \oplus P_{22}
\cr & \cr R_{4k+1} \cong & P_{11}^{2k+1} \oplus P_{12} \oplus P_{21}
\cr &\cr R_{4k+2} \cong & P_{11}^{2k+1} \oplus P_{12} \oplus P_{21}
\cr &\cr R_{4k+3} \cong & P_{11}^{2k+2} \oplus P_{22}.
\end{array}
\end{equation}

\subsection{From one-sided to two-sided resolutions}
The  construction of the differential $\delta$ is based on the
construction of a minimal projective resolution of $A/ {\rm rad} A$
as a  right $A$-module as developed in
\cite{Green-Solberg-Zacharia}. These differentials in the
right module resolutions of $A/ \rad A$ will suggest how to explicitly
construct
 the differential of a minimal $A$-$A$-bimodule resolution of
$A$. Our approach follows that of \cite{Green-Snashall} (see also
\cite{Snashall-Taillefer}) and, for Koszul algebras, that of
\cite{Green-Hartman-Marcos-Solberg}.

We  recall the construction of projective
resolutions as developed in \cite{Green-Solberg-Zacharia}, see also
\cite{Green-Snashall}.
Note that our labelling is different,  we use the
index to label the degree, and the exponent to count elements in a
fixed degree.  Let
$\Lambda$ be a finite dimensional algebra of the form
$\Lambda=KQ/I$. In \cite{Green-Solberg-Zacharia} the authors then
define sets  of paths $g_n$ in $KQ$ for $n \geq 0$ where $g_0$ is
given by the trivial paths in $KQ$, $g_1$ corresponds to the arrows
of $Q$ and $g_2$ is a minimal set of uniform generators of $I$. For
$n \geq3$, $g_n$ is constructed recursively such that the elements
in $g_n$ are uniform and satisfy the property: for $x \in g_n$,  $x=
\sum_{p \in g_{n-1}} p \lambda_p = \sum_{q \in g_{n-2}} q \mu_q$ for
unique elements $\lambda_p, \mu_q \in KQ$ and $\mu_q \in I$, having
special properties related to the minimal projective
right resolution of $A/{\rm rad} A$. If so then
there is a minimal projective $\Lambda$-module resolution
$$P_\bullet: \;\;\;\;\;\;\;  \ldots \longrightarrow P_n \stackrel{d_n}{\longrightarrow} P_{n-1}
\stackrel{d_{n-1}}{\longrightarrow} \ldots \longrightarrow P_1
\stackrel{d_1}{\longrightarrow} P_0 \longrightarrow \Lambda/\rm{rad}
\Lambda \longrightarrow 0$$ such that
\begin{enumerate}
\item $P_n = \bigoplus_{x \in g_n } {\bf t}(x) \Lambda$ for all $n \geq
0$
\item the map $d_n: P_n  \rightarrow P_{n-1} $ is given by
$${\bf t}(x) a \mapsto \sum_{p \in g_{n-1}} \lambda_p {\bf t}(x) a.
$$
In matrix terms the $(p,x)$-entry is $\lambda_p$.
\end{enumerate}
The collection of elements of $g_n$ is written as $g_n=\{g^i_n\} $.

\bigskip
In \cite{Green-Snashall} the sets $g_n$ are then used to give an
explicit construction of the first three differentials of a minimal
projective $\Lambda$-$\Lambda$-bimodule resolution of $\Lambda$.
Following their approach we first construct the sets $g_n$ for the
minimal right $A$-module resolution of $A/\rm{rad} A$ and then
proceed to construct an $A$-$A$-bimodule resolution of $A$.

\subsection{A minimal projective right $A$-module resolution of $A/{\rm rad}
A$}

In this paragraph we construct a minimal projective right $A$-module
resolution of $A/{\rm rad} A$. To facilitate the notations we
slightly change the notation of \cite{Green-Solberg-Zacharia} and
define the sets defining the differentials as follows:
\begin{Definition}\label{LH}
 For the algebra $A$ we set $g_0^1 =e_1$ and
$f^{22}_0 = e_2$ and we define recursively, for all $i= 0,1,2,3$ and
all $k \geq 0$ and for all $n \geq 4$,
$$
\begin{array}{lll}
 g^r_n &=& \left\{ \begin{array}{ll}
g^1_{n-1}\e - g^2_{n-1} \alpha \beta & \;\;\;\; \mbox{if $r=1$} \\
g^{r-1}_{n-1} \e - g^{r+1}_{n-1} \alpha \beta & \;\;\;\; \mbox{if
$r$ odd,  $1 < r \leq 2k-1$} \\ g^{r+1}_{n-1} \e + g^{r-1}_{n-1}
\alpha \beta & \;\;\;\; \mbox{if
$r$ even, $1 < r \leq 2k-1$ or $n=4k+2, r=2k$} \\
\end{array} \right. \\ \\ \\
 g^{2k+1}_{4k+i} &=& \left\{ \begin{array}{ll}
 g^{2k+1}_{4k+i-1}\e  & \;\;\;\; \mbox{if $i=1$ and $k=0$}\\
g^{2k}_{4k+i-1}\e  & \;\;\;\; \mbox{if $i=0,1$ and $k>0$} \\
g^{2k+1}_{4k+i-1} \e - f^{12}_{4k+i-1} \beta \alpha \beta  &
\;\;\;\;
\mbox{if $i=2,3$ and $k=0$} \\
g^{2k}_{4k+i-1} \e - f^{12}_{4k+i-1} \beta \alpha \beta  & \;\;\;\;
\mbox{if $i=2,3$ and $k >0$} \\
\end{array} \right. \\ \\ \\
 g^{2k}_{4k}&=&  g^{2k-1}_{4k-1} \alpha \beta   \;\;\;\; \mbox{if $k >0$} \\ \\ \\
 g^{2k+2}_{4k+3}&=&  g^{2k+1}_{4k+2} \alpha \beta \;\;\;\; \mbox{if $k >0$} \\ \\ \\
 f^{12}_{4k+i} &=&  \begin{array}{ll}
 g^{2k+1}_{4k+i-1} \alpha & \;\;\;\; \mbox{if
$i =1,2$}
\end{array} \\ \\ \\

f^{21}_{4k+i} &=&  \left\{ \begin{array}{ll}
 f^{22}_{4k} \beta & \;\;\;\; \mbox{if
$i =1$} \\ f^{21}_{4k+1} \e & \;\;\;\; \mbox{if $i =2$}
\end{array}  \right. \\ \\ \\

f^{22}_{4k+i} &=&   \left\{ \begin{array}{ll}
 f^{22}_{4k-1} (\beta \alpha)^2  & \;\;\;\; \mbox{if
$i =0$} \\ f^{21}_{4k+2} \alpha & \;\;\;\; \mbox{if $i =3$}
\end{array} \right. \\ \\ \\
\end{array}
$$

%Throughout this paper we set the convention that $g_{2k+m}^{4k+i} =
%0$ if $m>1$ and $i = 0,1,2$ or $m>2$ and $i=3$ and $g_0^n = 0 $ for
%all $n$.
We define $\cal{G}_n$ to be the set containing all elements of the
form $g_n^*$ and $f^{*}_n$.
\end{Definition}
{\bf Remark.}We note that $g_1^1 = \e$, $f^{12}_1 = \alpha$ and
$f^{21}_1 = \beta$, hence the set corresponding to $\cal{G}_1$ is
the set of arrows of ${\mathcal Q}$. Also $g^1_2 = \e^2 - \alpha
\beta \alpha \beta $, $f^{12}_{2} = \e \alpha$ and $f^{21}_2 = \beta
\e$ so that the set corresponding to  $\cal{G}_2$ is a minimal set
of uniform elements generating $I$.

Furthermore, we have chosen our notation such that it indicates from
which vertex to which other vertex the elements go, for example:
$g^r_{n-1}, g^r_{n}, g^r_{n+1}, {\it etc.}$ are all a sum of paths
from vertex one to itself, $f^{12}_{n-1}, f^{12}_{n}, f^{12}_{n+1}$,
{\it etc.}  go from vertex one to vertex two, and so on.

\bigskip

Let $P_n$ and $d_n$ be as defined in section 2.1.
It then
follows directly from \cite{Green-Solberg-Zacharia} that
\begin{Proposition}
The complex  $(P_\bullet, d)$ is a minimal projective right
$A$-module resolution of $A/ {\rm rad} A$.
\end{Proposition}

\subsection{A bimodule resolution} We will below
define differentials, and it is straightforward to check that they
define a complex. In order to show that it is exact, that is it
defines a projective $A$-$A$-bimodule resolution of $A$, we will use
the sets $\cal{G}_n$ defined the previous paragraph (in fact they
were crucial to find the differentials).

% We note that we can write $g_1^3 = g_1^2 \e
%- f_{12}^{2} \beta \alpha \beta = \e g_1^{2} - \beta \alpha \beta
%f_{12}^{2}$ and $g_1^{4k+i} = g_1 \e - g_2 \alpha \beta = \e g_1 -
%\alpha \beta g_3$. Continuing in this way we will be able to define
%the differential of a projective $A$-$A$-bimodule resolution of $A$.

We need to introduce some notation: Let $P_{ij} = Ae_i\otimes e_jA$,
we write the elements
$e_i \otimes e_j$ of $P_{ij}$  as follows:
$$ \begin{array}{lllll}
\aa & := &  e_1 \otimes e_1 &\in& P_{11} \\
\bb & := &  e_1 \otimes e_2 &\in& P_{12}\\
\cc & := &  e_2 \otimes e_1 &\in& P_{21} \\
\dd & := &  e_2 \otimes e_2 &\in& P_{22}. \\
\end{array}
$$
The terms of $R_\bullet$ in each degree contain multiple copies of
$P_{11}$ and always at most one copy of $P_{12}$, $P_{21}$ and
$P_{22}$ (see (\ref{terms})).
%In order to define the differential $\delta$ of $(R_\bullet, \delta)$
We need to distinguish the different copies of $P_{11}$ in $R_n$,
for $n \geq 0$, and we do this by adding superscripts. As before, the
subscripts determine in which degree of the complex an element lies.
For example,  $\aa_{i}^r = e_1 \otimes e_1 $ lies in the $r$th copy
of $P_{11}$ as a summand of $R_i$.

%, $\aa_{i+1}^{r+1} = e_1 \otimes
%e_1$ lies in the $r+1$th copy of $P_{11}$ as a summand of $R_{i+1}$
%and $\bb_{4k+1} = e_1 \otimes e_2$ (resp. $\cc_{4k+1} = e_2 \otimes
%e_1$) lies in the copy of $P_{12}$ (resp. $P_{21}$) as a summand of
%$R_{4k+1}$.

%$R_{4k} = \bigoplus_{r = 1}^{2k+1} P_{11}^r \oplus P_{22}$ and
%$\aa_{4k}^r$ denotes the element $e_1 \otimes e_1$ in the $r$th
%component $P_{11}^r$ of $R_{4k}$.

\begin{Definition}\label{differentials}
The map $\delta_0: R_0 \rightarrow A$ is given by multiplication.
With  the above notation, for $n \geq 1$, we define the differential
$\delta_n: R_n \rightarrow R_{n-1}$ recursively as follows.
$$ \aa^r_n \mapsto \left\{
\begin{array}{ll}
\aa^1_{n-1} \e - \e \aa^1_{n-1} & n=1, r=1 \\
\aa^1_{n-1} \e - \bb_{n-1} \beta \alpha \beta + \e \aa^1_{n-1} -
\alpha \beta \alpha \cc_{n-1} - \alpha \beta \bb_{n-1} \beta -
\alpha \cc_{n-1} \alpha \beta &
n=2, r=1 \\
\aa^1_{n-1} \e - \bb_{n-1} \beta \alpha \beta - \e \aa^1_{n-1} +
\alpha \beta \alpha
\cc_{n-1} & n=3, r=1 \\
\aa^1_{n-1} \alpha \beta - \alpha \beta \aa_{n-1}^1 - \e \bb_{n-1}
\beta + \alpha
\cc_{n-1} \e & n=3, r=2 \\
\end{array} \right.
$$

and, for all $n = 4k+i$ with $i = 0,1,2,3$ and $k \geq 1$:

If $i=0, 2$
$$ \aa^r_{4k+i} \mapsto \left\{
\begin{array}{ll}
\aa^1 \epsilon -\aa^2 \alpha \beta + \epsilon \aa^1 -\alpha \beta \aa^2 & \mbox{if $r=1$} \\
 \aa^{r+1} \e  + \aa^{r-1} \alpha \beta  + \e \aa^r - \alpha \beta \aa^{r+2} & \mbox{if $ 1 < r < 2k$, $r$ even} \\
 \aa^{r-1} \e  - \aa^{r+1} \alpha \beta  + \e \aa^r + \alpha \beta \aa^{r-2} & \mbox{if $ 1 < r < 2k$, $r$ odd} \\
\aa^{r-1} \alpha \beta  + \e \aa^r + (-1)^k \alpha \beta \alpha \dd \beta & \mbox{if $r=2k$ and $i=0$} \\
 \aa^{r-1} \e    + \alpha \beta \aa^{r-2} + (-1)^k  \alpha  \dd \beta \alpha \beta& \mbox{if $r=2k+1$ and $i=0$} \\
\aa^{r+1}\e + \aa^{r-1} \alpha \beta  + \e \aa^r + (-1)^{k+1} \alpha \beta \alpha \cc -\alpha \beta \bb \beta & \mbox{if $r=2k$ and $i=2$} \\
 \aa^{r-1} \e  + \e \aa^r  + \alpha \beta \aa^{r-2} - \bb \beta \alpha \beta+ (-1)^{k+1}  \alpha  \cc  \alpha \beta& \mbox{if $r=2k+1$ and $i=2$}
\end{array}
\right.
$$
If $i=1, 3$

$$ \aa^r_{4k+i} \mapsto \left\{
\begin{array}{ll}
\aa^{1} \e  - \aa^{2} \alpha \beta  - \e \aa^1 + \alpha \beta \aa^{3} & \mbox{if $r=1$ } \\
\aa^{3} \e  + \aa^{1} \alpha \beta  - \e \aa^2 - \alpha \beta \aa^{1} & \mbox{if $r=2$ } \\
\aa^{r-1} \e  - \aa^{r+1} \alpha \beta  - \e \aa^r + \alpha \beta \aa^{r+2} & \mbox{if $ 1 < r \leq 2k$ and $r$ is odd} \\
\aa^{r+1} \e  + \aa^{r-1} \alpha \beta  - \e \aa^r - \alpha \beta \aa^{r-2} & \mbox{if $ 1 < r \leq 2k$ and $r$ is even} \\
\aa^{r-1} \e   - \e \aa^{r}  & \mbox{if $r=2k+1$ and $i=1$} \\
\aa^{r-1} \e   - \bb \beta \alpha \beta - \e \aa^{r} + (-1)^k  \alpha \beta \alpha \cc             & \mbox{if $r=2k+1$ and $i=3$} \\
\aa^{r-1} \alpha \beta   -  \alpha \beta \aa^{r-2} + (-1)^k  \alpha \cc \e - \e  \bb \beta          & \mbox{if $r=2k+2$ and $i=3$} \\
\end{array} \right. $$
where on the right hand side we omit the subscript with the
understanding that these are all elements in $R_{4k+i-1}$. Finally,
for $n= 4k+i$ with $i = 0,1,2,3$ and $k \geq 0$, we set
$$
\begin{array}{lll}
\dd_{4k} & \mapsto &  \dd_{4k-1} (\beta \alpha)^2 + (\beta \alpha)^2
\dd_{4k-1}
+ (-1)^k \beta \aa_{4k-1}^{2k-1} \alpha \\ \\
\bb_{4k+1} & \mapsto & \aa^{2k+1}_{4k} \alpha +(-1)^{k+1} \alpha \dd_{4k} \\ \\
\cc_{4k+1} & \mapsto & \dd_{4k} \beta +(-1)^{k+1} \beta \aa^{2k}_{4k} \\ \\
\bb_{4k+2} &  \mapsto & \aa^{2k+1}_{4k+1} \alpha + \e \bb_{4k+1} \\ \\
\cc_{4k+2} &\mapsto& \cc_{4k+1} \e +(-1)^{k} \beta \aa^{2k+1}_{4k+1} \\ \\
\dd_{4k+3} &\mapsto& \cc_{4k+2} \alpha + (-1)^{k+1} \beta \bb_{4k+2}. \\ \\
\end{array}$$

\end{Definition}

Recall from Definition~\ref{LH} that the set $\cal{G}_n$ consists of
all the elements of the from $g_n^*$ and $f_n^*$. Furthermore, as
stated in the remark following definition~\ref{LH} we know the
origin and target of each element of the set $\cal{G}_n$. A simple
counting argument now shows that if we set  $R_n = \bigoplus_{x \in
{\mathcal G}_n} A{\bf o}(x) \otimes {\bf t}(x) A$ we obtain exactly
the expressions (2) of  section 2. Furthermore, we have that,

\begin{Theorem}
The complex $(R_\bullet, \delta)$ is a minimal projective resolution
of $A$ as an $A$-$A$-bimodule.
\end{Theorem}
{\it Proof:} It is straightforward to verify that $\Ima \ \delta_n
\subseteq \Ker \ \delta_{n-1}$ and therefore  $(R_\bullet, \delta)$ is
a complex.

In order to show that the complex is exact we first note that there
is an isomorphism of complexes of right $A$-modules $(A/ \rad A
\otimes_A R_\bullet, id \otimes \delta) \simeq (P_\bullet, d)$. This
holds since - as it is easily verified -  for all $n$, we have $A/ \rad
A \otimes_A R_n \simeq P_n$ and furthermore, the diagram

\begin{equation}\label{equivcomplexes}
\xymatrix{
 A/ \rad A \otimes_A R_n   \ar[d]^{\simeq} \ar[r]^{id \otimes \delta_n} &A/ \rad A \otimes_A R_{n-1}\ar[d]^{\simeq}\\
P_n \ar[r]^{d_n} &P_{n-1}}
\end{equation}

commutes for all $n \geq 1$. Now we apply the same arguments as in
\cite[Proposition 2.8]{Green-Snashall} and see also \cite[Theorem
1.6]{Snashall-Taillefer}. Suppose that $\Ker \  \delta_{n-1}
\not\subseteq \Ima \ \delta_n$ for some $n \geq 1$. Then there exists
a simple $A$-$A$-bimodule $S \otimes T$ (where $S$ is a
simple left
$A$-module and $T$ is  a simple right $A$-module) and there exists a
non-zero map $f: \Ker \ \delta_{n-1} \rightarrow S \otimes T$. We have
seen (Proposition 2.2 and~(\ref{equivcomplexes})) that the complex
$(A/ \rad A \otimes_A R_\bullet, id \otimes \delta)$ is a minimal
projective resolution of $A /\rad A$ as a right $A$-module.
Therefore the following isomorphisms of right $A$-modules hold $
A/\rad A \otimes_A \Ima \  \delta_n \simeq \Ima (id \otimes \delta_n)
\simeq \Ker (id \otimes \delta_{n-1}) \simeq A /\rad A \otimes_A
\Ker \delta_{n-1}$ and we obtain a non-zero map by taking  the composition

$$  A/ \rad A \otimes_A R_n  \stackrel{id \otimes
\delta_n}{\longrightarrow} A/\rad A \otimes_A \Ima \ \delta_n
\stackrel{\simeq}{\longrightarrow} A /\rad A \otimes_A \Ker
\ \delta_{n-1} \stackrel{id \otimes f}{\longrightarrow}S \otimes T.$$

But this is the same as the functor $A/\rad A \otimes_A - $ applied
to the sequence of maps $R_n \stackrel{\delta_n}{\rightarrow} \Ker
\ \delta_{n-1} \stackrel{f}{\rightarrow} S \otimes T$. However,  the
map $f \circ \delta_n$ is zero. This gives a contradiction and
therefore the complex $(R_\bullet, \delta_n)$ is exact.

Since we calculated the terms of $R_\bullet$ following
\cite{Happel}, the minimality follows.
 $\Box$

\section{Dimensions of the Hochschild cohomology groups}

The Hochschild cohomology of $A$ is defined by
$\HH^*(A)=\Ext^*_{A^e}(A,A)$. From now we denote by $\Omega$ the syzygy
operator for $A^e$-modules,
it  associates to an $A^e$-module $M$
the $A^e$-module given by the kernel of the map $P(M) \rightarrow M$
where $P(M)$ is the projective cover of $M$ and iterating this, defines
$\Omega^n$. Then for all $n \geq 0$, there is an exact
sequence

\begin{equation}\label{exactsequence}
0 \rightarrow \Hom_{A^e}(\Omega^n(A),A) \rightarrow
\Hom_{A^e}(R_n,A) \rightarrow  \Hom_{A^e}(\Omega^{n+1}(A),A)
\rightarrow \Ext^{n+1}_{A^e}(A,A) \rightarrow 0.
\end{equation}

and since we have an explicit description of the terms $R_n$ of the
$A$-$A$-bimodule resolution of $A$, it is straightforward to
calculate the dimensions of the complex $\Hom_{A^e}(R^\bullet, A)$.
These dimensions are for, $k \geq 0$,

$$
\begin{array}{lcl}
\dim \Hom_{A^e} (R_{4k},A ) &=& 8k+7 \\
\dim \Hom_{A^e} (R_{4k+1},A ) &=& 8k+8 \\
\dim \Hom_{A^e} (R_{4k+2},A ) &=& 8k+8 \\
\dim \Hom_{A^e} (R_{4k+3},A ) &=& 8k+11. \\
\end{array}
$$

Furthermore,  $\HH^0(A)$ is given by the center $Z(A)$ of $A$ and this
has vector space basis
$$ \{e_1 + e_2, \varepsilon, \alpha \beta + \beta \alpha, \varepsilon^2, \beta\alpha\beta\alpha\}.$$

By (\ref{exactsequence}), in order to calculate the dimensions of
$\HH^{n}(A)$ for all $n \geq 0$ we need to calculate the dimensions
of $\Hom_{A^e}(\Omega^n(A),A)$ for all $n \geq 1$. Using the
description

$$ \Hom_{A^e}(\Omega^n(A),A) = \{ \gamma \in \Hom_{A^e}(R_n,A) | \;\; \gamma(\Omega(A)^{n+1}) = 0 \} $$

we will calculate a basis for $\Hom_{A^e}(\Omega^n(A),A)$.
 Since $R_\bullet$ is exact the images of the differential yield
 the generators of $\Omega^n(A)$ in each degree $n$. Thus, for $k \geq 0$, we have

 $$ \begin{array}{lll}
\Omega^{4k}(A) &=& \langle \delta_{4k}(\aa^1_{4k}), \ldots, \delta_{4k}(\aa^{2k+1}_{4k}),  \delta_{4k}(\dd_{4k}) \rangle \\ \\
\Omega^{4k+1}(A) &=& \langle \delta_{4k+1}(\aa^1_{4k+1}), \ldots,
\delta_{4k+1}(\aa^{2k+1}_{4k+1}),
\delta_{4k+1}(\bb_{4k+1}),  \delta_{4k+1}(\cc_{4k+1})\rangle \\ \\
\Omega^{4k+2}(A) &=& \langle \delta_{4k+2}(\aa^1_{4k+2}), \ldots,
\delta_{4k+2}(\aa^{2k+1}_{4k+2}),
\delta_{4k+2}(\bb_{4k+2}),  \delta_{4k+2}(\cc_{4k+2}) \rangle\\ \\
\Omega^{4k+3}(A) &=& \langle \delta_{4k+3}(\aa^1_{4k+3}), \ldots, \delta_{4k+3}(\aa^{2k+2}_{4k+3}),  \delta_{4k+3}(\dd_{4k+3}) \rangle.\\ \\
\end{array}$$

 For each of the four cases above, we will obtain a system of equations whose solution space gives a basis of
 $\Hom_{A^e}(\Omega^n(A),A)$.

\subsection{Characteristic $\neq 2$}

Suppose that the characteristic of $K$ is $\neq 2$. Then

\begin{Lemma}\label{dimOmega}
We have the following dimensions, for $k \geq 0$:
$$\begin{array}{lll}
\dim \Hom_{A^e}(\Omega^{4k+i}(A),A)  &=& \left\{ \begin{array}{ll}
5k+5 &  \;
\mbox{if}  \;\;\; i=0\\
 5k+5 &  \; \mbox{if}  \;\;\; i=1 \\
 5k+6 &  \; \mbox{if}  \;\;\;  i=2 \\
5k+6 &  \; \mbox{if}  \;\;\; i=3 \\
\end{array} \right.
\end{array}$$
\end{Lemma}

\bigskip
{\it Proof:}
 Let $k \geq 1$, we start by
calculating the dimension of $\Hom_{A^e}(\Omega^{4k}(A),A)$. Let
$\gamma: R_{4k} \rightarrow A$ be a most general map defined by
$\gamma (\aa_{4k}^j) = a_1^j e_1 + b_1^j \e + c_1^j \alpha
\beta + d_1^j \e^2 $ for $1 \leq j \leq 2k +1$ and $\gamma
(\dd_{4k}) = a_2 e_2 + b_2 \beta \alpha + c_2 \beta \alpha \beta
\alpha$ with $a_1^j, b_1^j, \ldots, a_2, b_2 \in k$. With this
notation we have that $\gamma(x) = 0$ for all generators $x$ of
$\Omega^{4k+1}(A)$ if and only if for all even $n$, such that $2
\leq n \leq 2k$,

$$\begin{array}{rllcrllcrll}
a_1^n &=& a_1^{n+1}& \;\;\;\;\;\; & b_1^n &=& b_1^{n+1}&
\;\;\;\;\;\; & c_1^n &=& c_1^{n+1}
 \\
\end{array}$$
and
$$(-1)^{k}a_2 = a_1^{2k+1}$$ and $$(-1)^{k}b_2 =c_1^{2k+1}.$$

Thus the dimension of $\Hom_{A^e}(\Omega^{4k}(A),A)$ is $5k+5$.

\bigskip

Similarly for the dimension of $\Hom_{A^e}(\Omega^{4k+1}(A),A)$, let
$\gamma: R_{4k+1} \rightarrow A$ be the map defined by
$\gamma (\aa_{4k+1}^i) = a_1^i e_1 + b_1^i \e + c_1^i \alpha
\beta + d_1^i \e^2 $ for $1 \leq i \leq 2k +1$ and $\gamma
(\bb_{4k+1}) = a_{12} \alpha + b_{12} \alpha \beta \alpha$ and
$\gamma(\cc_{4k+1}) = a_{21} \beta + b_{21} \beta \alpha
\beta $ with $a_1^i, b_1^i, \ldots, a_{12}, b_{21}, \ldots \in k$.
Then for all generators $x$ of $\Omega^{4k+2}(A)$, we have that
$\gamma(x) = 0$ if and only if the following equations hold

\begin{equation*} 2 a_1^1 = 2 a_1^2 = 0 \tag{*} \end{equation*}
$$ a_1^3 = a_1^4 = \ldots =a_1^{2k+1} = a_1^{2k+1} =0$$
\begin{equation*} 2 ( b_1^{1} - c_1^{2})=0 \tag{*} \end{equation*}
$$b_1^{2k} + b_1^{2k+1} + c_1^{2k-1} - a_{12} + (-1)^{k+1} a_{21} = 0$$
$$ \begin{array}{rll}
c_1^{2k+1}&=&0 \\
b_1^{n} +   b_1^{n+1} +   c_1^{n-1} -  c_1^{n+2}  &=& 0 \\
\end{array}$$

 for $n \in \{2,4, \ldots, 2k-2\}$.

Thus the dimension of $\Hom_{A^e}(\Omega^{4k+1}(A),A)$ is $5k+5$.

In the same way we obtain the following sets of equations for
$\Hom_{A^e}(\Omega^{4k+2}(A),A)$ for all even $n$, such that $2 \leq
n \leq 2k$,

$$\begin{array}{rllcrllcrll}
a_1^n &=& a_1^{n+1}& \;\;\;\;\;\; & b_1^n &=& b_1^{n+1}&
\;\;\;\;\;\; & c_1^n &=& c_1^{n+1}
 \\
\end{array}$$

and

$$a_{12} = (-1)^{k}a_{21}  \;\; \mbox{and} \;\;  b_{12} = (-1)^{k}b_{21} $$

giving $\dim \Hom_{A^e}(\Omega^{4k+2}(A),A) =5k+6$.

Finally for $\Hom_{A^e}(\Omega^{4k+3}(A),A)$ we obtain

\begin{equation*} 2 a_1^1 = 2 a_1^2 = 0 \tag{*} \end{equation*}
$$a_1^3 = a_1^4 = \ldots =a_1^{2k+2} = 0$$
\begin{equation*} 2 ( b_1^{1} - c_1^{2})=0 \tag{*} \end{equation*}
$$ \begin{array}{lll}
b_1^{n} +   b_1^{n+1} +   c_1^{n-1} -  c_1^{n+2}  &=& 0 \\
\end{array}$$

 for $n \in \{2,4, \ldots, 2k\}$ and

$$b_1^{2k+2} = (-1)^{k} a_2 $$ and

\begin{equation*} c_1^{2k+2} = (-1)^{k}  2a_2   \tag{*} \end{equation*}

such that $\dim \Hom_{A^e}(\Omega^{4k+3}(A),A) =5k+6$.

\bigskip

From Definition~\ref{differentials} we see that, for $n \leq 3$, the
pattern for the expressions of the differential is not regular.
Therefore in these small cases we have to calculate the dimension of
$\Hom_{A^e}(\Omega^{n}(A),A)$ by hand. These calculations show that
$\dim \Hom_{A^e}(\Omega^{1}(A),A) =5$, $\dim
\Hom_{A^e}(\Omega^{2}(A),A) =6$ and $\dim
\Hom_{A^e}(\Omega^{3}(A),A) =6$.
 $\Box$

\bigskip

The dimensions of the Hochschild cohomology spaces of $A$
follow directly from above.
\bigskip
\begin{Theorem}
Assume char $K \neq 2$. The dimensions of the Hochschild cohomology
spaces of $A$
%in small degrees are given by
%\begin{center}
%\begin{tabular}{l||l|l|l|l|l|l|l|l|l|l|l}
%i &0 & 1 & 2 &3 & 4& 5& 6& 7& 8& 9& 10  \\  &&&&&&&&&&&\\
%\hline   &&&&&&&&&&&\\
%$ \dim \HH^i(A)$  &4 &2 &3 & 5&  6& 1& 5& 8& 9& 7& 7                     \\
%\end{tabular}
%\end{center}
%\bigskip
%Starting from $\HH^{11}(A)$ the Hochschild cohomology
%are periodic of period 4 and the dimensions
are as follows
$$ \dim \HH^{4k+i}(A) = \left\{
\begin{array}{ll}
2k+3 & \mbox{$i =0, k>0$ } \\
2k+3 & \mbox{$i =1, k\geq 0 $} \\
2k+3 & \mbox{$i =2, k\geq 0$ } \\
2k+4 & \mbox{$i =3, k \geq 0$ } \\
\end{array}
\right. $$ and the set $\{e_1+e_2, \e, \alpha \beta + \beta \alpha,
\e^2, \beta \alpha \beta \alpha \}$ is a $K$-basis of  $\HH^0(A)$.
\end{Theorem}

\subsection{Characteristic of  $K = 2$}

Suppose the characteristic of $K =2$. Then

\begin{Lemma}
We have the following dimensions, for $k \geq 0$:

$$\begin{array}{lll}
\dim \Hom_{A^e}(\Omega^{4k+i}(A),A)  &=& \left\{ \begin{array}{ll}
5k+5 &  \; \mbox{if}  \;\;\; i=0\\
 5k+8 &  \; \mbox{if}  \;\;\; i=1, k \neq 0 \\
 5k+6 &  \; \mbox{if}  \;\;\;  i=2 \\
5k+9 &  \; \mbox{if}  \;\;\; i=3 \\
\end{array} \right.
\end{array}$$
and $\dim \Hom_{A^e}(\Omega^{1}(A),A) =6$.
\end{Lemma}

\bigskip
{\it Proof:} \ The proof is analogous to the proof of
Lemma~\ref{dimOmega} except in the rows marked with (*)
we either gain or loose a
degree of freedom. This only occurs in the cases $i=1$ and $i=3$.
Furthermore, as before we have to calculate the cases of small
degree separately by hand. \hfill $\Box$

\bigskip

As before, the dimensions of the Hochschild cohomology spaces  of
$A$ directly follow.

\begin{Theorem}
Assume char $K = 2$. The dimensions of the Hochschild cohomology
spaces of $A$
%in small degrees are given by
%\begin{center}
%\begin{tabular}{l||l|l|l|l|l|l|l|l|l|l|l}
%i &0 & 1 & 2 &3 & 4& 5& 6& 7& 8& 9& 10  \\  &&&&&&&&&&&\\
%\hline   &&&&&&&&&&&\\
%$ \dim \HH^i(A)$  &4 &2 &3 & 5&  6& 1& 5& 8& 9& 7& 7                     \\
%\end{tabular}
%\end{center}
%\bigskip
%Starting from $\HH^{11}(A)$ the Hochschild cohomology
%are periodic of period 4 and the dimensions
are as follows
$$ \dim \HH^{4k+i}(A) = \left\{
\begin{array}{ll}
2k+6 & \mbox{$i =0, k>0$ } \\
2k+6 & \mbox{$i =1, k > 1 $} \\
2k+6 & \mbox{$i =2, k > 0$ } \\
2k+7 & \mbox{$i =3, k > 0$ } \\
\end{array}
\right. $$ and $ \dim \HH^{1}(A) = 4$. The set $\{e_1+e_2, \e,
\alpha \beta + \beta \alpha, \e^2, \beta \alpha \beta \alpha \}$ is
a $K$-basis of $\HH^0(A)$.
\end{Theorem}

{\bf Acknowledgements.} The authors would like to thank Nicole
Snashall and the referee for their very helpful suggestions and
comments.

\end{document}